\documentclass{article}
\usepackage{amsmath,amsthm,amsopn,amstext,amscd,amsfonts,amssymb,verbatim}

\usepackage{color}
\usepackage{natbib}
\usepackage{graphicx}
\def\tr{\mathop{\rm tr}\nolimits}

\def\diag{\mathop{\rm diag}\nolimits}

\pdfoutput=1

\renewenvironment{abstract}
                 {\vspace{6pt}
                  \begin{center}
                  \begin{minipage}{5in}
                  \centerline{\textbf{Abstract}}
                  \noindent\ignorespaces
                 }
                 {\end{minipage}\end{center}}
\newtheorem{thm}{\textbf{Theorem}}[section]

\theoremstyle{definition}
\newtheorem{defn}{\textbf{Definition}}[section]

\title{\Large \textbf{Polynomial Eulerian shape distributions}}
\author{
  \textbf{Francisco J. Caro-Lopera} \thanks{Corresponding author\newline
    }\\
  {\normalsize Departament of Basic Sciences} \\
  {\normalsize Universidad de Medell\'{\i}n} \\
  {\normalsize Medell\'{\i}n, Colombia} \\
  {\normalsize E-mail: fjcaro@udem.edu.co} \\[2ex]
  \textbf{Jos\'e A. D\'{\i}az-Garc\'{\i}a} \\
  {\normalsize Universidad Aut\'onoma Agraria Antonio Narro} \\
  {\normalsize Calzada Antonio Narro 1923, Col. Buenavista} \\
  {\normalsize 25315 Saltillo, Coahuila, M\'exico}\\
  {\normalsize E-mail: jadiaz@uaaan.mx}\\
}
\date{}
\begin{document}
\maketitle

\begin{abstract}
In this paper a new approach is derived in the context of shape theory. The implemented
methodology is motivated in an open problem proposed in \citet{GM93} about the construction of
certain shape density involving Euler hypergeometric functions of matrix arguments.

The associated distribution is obtained by establishing a connection between the required shape
invariants and a known result on squared canonical correlations available since 1963; as usual
in statistical shape theory and the addressed result, the densities are expressed in terms of
infinite series of zonal polynomials which involves considerable difficulties in inference.
Then the work proceeds to solve analytically the problem of computation by using the Eulerian
matrix relation of two matrix argument for deriving the corresponding polynomial distribution
in certain parametric space which allows to perform exact inference based on exact
distributions characterized for polynomials of very low degree. A methodology for comparing
correlation shape structure is proposed and applied in handwritten differentiation.
\end{abstract}

\section{Introduction}\label{sec:introduction}
Invariant statistics has emerged as one of the most powerful tools in statistics to study
matrix variate problems, in particular the statistical shape analysis has been rolled a number
of applications in almost all the natural and exact sciences; their applications in machine
vision and images analysis are notable. In general, shape theory deals with problems of
differentiation of ``objects" which can be summarized by anatomical or mathematical marks
called landmarks, under the assumption that those points can be placed (in presence of
randomness) in any similar objects and absorbs all the geometrical information of the ``shape"
of the object. The invariant theory enters in the analysis, when the comparisons, matching,
and/or inference are performed in special quotient spaces rather than the usual full noised
Euclidean space. Then according to the symmetries involved in the study defined and supervised
by the expert in the experiment, a number of different shapes can be defined. Three classical
shapes are well studied in literature, the Euclidean or similarity shape, the affine or
configuration, and the projective shape; then the associated shape is all geometrical
information which remains under removing location, rotation, scaling, reflection, uniform shear
and/or projection, etc., see for example, \citet{DM98} and the references therein.

Under deterministic assumptions, shape theory counts a number of approaches, a source of their
main aspects are given for example in  \citet{Kendalletal99}, their references; as a
consequence of this geometrical understanding of shape, the procrustes theory, including
projective shape analysis and its relation to affine shape analysis, have been studied
profusely in applications; some works in these directions are due to \citet{Mardiaetal05},
\citet{Patran03}, \citet{Maietal11}, \citet{DomokosKato2010}, \citet{ET04}, \citet{GlasMar01},
\citet{g:91}, \citet{GroTag09}, \citet{LinGar97}, \citet{Maretal96}, \citet{MarPatr05},
\citet{MokAbb02}, \citet{H:92}, \citet{K:04}, \citet{LinFang07}, among many others

Now, the statistical theory of shape deals with landmarks in presence of random (see for
example \citet{DM98}, \citet{Small96} and the references therein). A branch of the  statistical
shape analysis was connected by \citet{GM93}(corrected by \citet{dgr:03} and \citet{cdg:10}))
with the very old theory of matrix variate distribution based on A.T. James zonal polynomials
(see for example \citet{JAT64} and the cited papers). At that time the shape distributions
derived in \citet{GM93} and most of the classical results in matrix variate distributions
collected in the excellent book of \citet{MR1982} could not computed and force to use
asymptotic distributions. Only recent the zonal polynomial were computable in an efficient way
by \citet{KE06}, and a hope for computations of infinite series of zonal polynomials was a
possible task under certain restriction and subject to the problem of automatic detection of
convergence of the series.

Parallel to the problem of computation, the classical shape theory of \citet{GM93} based on the
Gaussian distributions, opened a field of new interesting theoretical and applied problems, the
so called generalized shape theory, which models shape under non Gaussian models, see for
example, \citet{Caro09b,cdg:10,cdg:14c},  \citet{dc:12a,dc:12b,dc:14,DC16} and \citet{cgb:13a},
among many others. Moreover the generalized shape analysis can be studied in the more general
setting of real normed division algebras, a unified approach valid for real, complex,
cuaternion and octonion landmark data matrices \citet{dc:10}.

As is can be check in the above references the Euclidean shape based on similarity
transformations, in all theirs forms, and the affine configuration have been studied
extensively, involving inference, families of distributions and so on; however, an open problem
seems to be unsolved from the times of \citet{GM93}.

Analysis of shape under Gaussian models were proposed by \citet{GM93} following a natural
induction way of deriving new shape distribution via Laplace-transform in order to obtain
hypergeometric series of matrix arguments by known integration over cones and manifolds
involving zonal polynomials. This methodology was though for obtaining successive distributions
which filters new more valuable information of the original landmark data.

The underlying sequence of new shape distributions appeared as follows: start with an original
Gaussian landmark data, then after removing traslation and rotation the so called size and
shape distribution is obtained, this function has the form of a ${}_{0}F_{1}(\cdot)$
hypergeometric function Bessel type (see \citet{MR1982}), a known series of zonal polynomials.
According to the inductive process of constructing hypergeometric series, the next appropriate
shape should be a confluent type  ${}_{1}F_{1}(\cdot)$ by removing the scaling information of
the size and shape, then the shape distribution is gotten. But as it can be seen in
\citet{GM93}, this shape has not the required form and no properties were available at that
time to explore the asymptotics or computation of the distribution (later
\citet{dc:12a,dc:12b,dc:14,DC16} generalized the result to every elliptical model and
transformations and performed inference with the exact densities). So, in order to obtain
confluent shape distribution, \citet{GM93} defined an affine transformation which lead directly
from the Helmertized matrix (the original Gaussian landmark matrix invariant under traslation)
to the required confluent distribution (later \citet{Caro09b,cdg:10,cdg:14c} generalized it to
the class of elliptical distribution and performed inference with the exact distribution, some
of them simple polynomials of low degree).  Finally, \citet{GM93} conjectured that a new class
of shape distributions could be constructed by using the induction way and claimed for an
hypergeometric type shape distribution ${}_{2}F_{1}(\cdot)$.

In this paper we study the addressed open problem by defining a new approach in the context of
statistical shape theory. Section \ref{sec:EulerianShape} focus in the definition of the new
shape invariant. As very usual in shape theory, the resulting density is a series of zonal
polynomials which presents strong problems for inference; it motivated Section
\ref{sec:polynomialEulerian} which solves the problem of inference by considering in some
simple parametric space an equivalent polynomial distribution of very low degree in the latent
roots of the matrix arguments of the zonal polynomials, then the exact formulae available for
these polynomials in two dimensions, a number of applications and methodologies for shape
comparisons are given in Section \ref{sec:applications}.

\section{Eulerian Shape}\label{sec:EulerianShape}
The aim of this approach consists of proposing an alternative way to compare populations in the
context of the shape theory. It is easy to check in the literature shape that Wilks' theorem is
the most common resource for testing means of two populations (a number of other tests can be
implemented, see \citet{DC08}), \cite{cgb:14a,cgb:14b}. However, the application of the
classical chi-squared distribution of that theorem demands a number of conditions on regularity
which are difficult to keep in most of the available landmark data.

The strength of relation between two variates has a similar power of explanation in shape
theory than the equality of means and both complement each other for completing a shape
analysis.

Consider two $N\times K$ matrices $\mathbf{X}^{**}$ and $\mathbf{Y}^{**}$, which summarize two
figures comprised in $N$ landmarks in $K$ dimensions, such that $\mathbf{X}^{**}\sim
\mathcal{N}(\boldsymbol{\mu}_{\mathbf{X}^{**}},\mathbf{I}_{N}\otimes
\boldsymbol{\Phi}_{\mathbf{X}^{**}})$ and $\mathbf{Y}^{**}\sim
\mathcal{N}(\boldsymbol{\mu}_{\mathbf{Y}^{**}},\mathbf{I}_{N}\otimes
\boldsymbol{\Phi}_{\mathbf{Y}^{**}})$, where $\boldsymbol{\Phi}_{\mathbf{X}^{**}}$ and
$\boldsymbol{\Phi}_{\mathbf{Y}^{**}}$ are positive definite $K\times K$ matrices.

Now, as usual in shape theory we remove some meaningless geometrical information as the
translation of both figures, by using contrast matrices, for example sub-Helmert arrays. In
this case we use a $(N-1)\times N$ matrix, $\mathbf{L}$, which has orthonormal rows to
$\textbf{1}=(1,...,1)'$; to construct the so called Helmertized $N-1\times K$ matrices, say,
$\mathbf{\mathbf{X}}^{*}$ and $\mathbf{Y}^{*}$, respectively, via
$\mathbf{LX}^{**}=\mathbf{\mathbf{X}}^{*}$ and $\mathbf{LY}^{**}=\mathbf{Y}^{*}$.

Under these transformations, the randomness of the transformed landmarks follow the laws:
$\mathbf{X}^{*}\sim N(\boldsymbol{\mu}_{\mathbf{X}^{*}},\mathbf{I}_{n}\otimes
\mathbf{\Sigma}_{22})$, $\mathbf{Y}^{*}\sim
\mathcal{N}(\boldsymbol{\mu}_{\mathbf{Y}^{*}},\mathbf{I}_{n}\otimes \mathbf{\Sigma}_{11})$,
where $n=N-1\times K$. The means are changed to
$\boldsymbol{\mu}_{\mathbf{X}^{*}}=\mathbf{L}\boldsymbol{\mu}_{\mathbf{X}^{**}}$ and
$\boldsymbol{\mu}_{\mathbf{Y}^{*}}=\mathbf{L}\boldsymbol{\mu}_{\mathbf{Y}^{**}}$, while the new
covariance matrices are conveniently denoted by $\mathbf{\Sigma}_{11}$ and
$\mathbf{\Sigma}_{22}$, where
$\mathbf{\Sigma}_{11}=\mathbf{L}\boldsymbol{\Phi}_{\mathbf{Y}^{**}}\mathbf{L}'$ and
$\mathbf{\Sigma}_{22}=\mathbf{L}\boldsymbol{\Phi}_{\mathbf{X}^{**}}\mathbf{L}'$.

We storage the above Helmertized matrices into an $N\times 2K$ partition array
$\mathbf{W}=[\mathbf{Y}^{*}|\mathbf{X}^{*}]$, which implies that $\mathbf{W}\sim
\mathcal{N}([\boldsymbol{\mu}_{\mathbf{Y}^{*}}|\boldsymbol{\mu}_{\mathbf{X}^{*}}],\mathbf{I}_{n}\otimes\mathbf{\Sigma})$,
where
$$
  \mathbf{\Sigma}=\left(\begin{array}{cc}\mathbf{\Sigma}_{11}&\mathbf{\Sigma}_{12}\\
  \mathbf{\Sigma}_{21}&\mathbf{\Sigma}_{22}\end{array}\right).
$$

Finally, set
$\mathbf{Z}=\mathbf{W}-[\boldsymbol{\mu}_{\mathbf{Y}^{*}}|\boldsymbol{\mu}_{\mathbf{X}^{*}}]$,
thus $\mathbf{Z}=[\mathbf{Y}|\mathbf{X}]\sim
\mathcal{N}(\mathbf{0},\mathbf{I}_{n}\otimes\mathbf{\Sigma})$, where we have defined
$\mathbf{Y}=\mathbf{Y}^{*}-\boldsymbol{\mu}_{\mathbf{Y}^{*}}$ and
$\mathbf{Y}=\mathbf{Y}^{*}-\boldsymbol{\mu}_{\mathbf{Y}^{*}}$.

Now,  we are interested in the measure of the relation between the two variates $\mathbf{X}$
and $\mathbf{Y}$, where $\mathbf{\Sigma}_{12}=\mathbf{\Sigma}_{21}'$ is the covariance matrix
between the components of $\mathbf{X}$ and the components of  $\mathbf{Y}$. In order to
establish that comparison we need to define the following supporting equivalent class of
shapes.

\begin{defn}\label{Def:eulerian}
Given two Helmertized $n\times s$ matrices $\mathbf{X}$ and $\mathbf{Y},$ it is said that they
have the same eulerian shape, if
\begin{enumerate}
    \item $\mathbf{HX}=[\mathbf{X}_{1}'|\,\, \mathbf{0}]'$, where $\mathbf{X}_{1}$ is a non-singular $s\times
    s$ matrix  and $\mathbf{H}\in O(r)$,
    \item $\mathbf{HY}=\mathbf{T}=[\mathbf{U}'|\mathbf{V}']'$, where $\mathbf{U}$ is $s\times s$, $\mathbf{V}$ is $(r-s)\times
s$ and $\mathbf{H}$ is defined according to (1), and,
    \item $\mathbf{U}$ and $\mathbf{V}$ satisfies $|(\mathbf{U}'\mathbf{U})(\mathbf{U}'\mathbf{U}+\mathbf{V}'\mathbf{V})^{-1}-\mathbf{I}_{s}|=0$.
\end{enumerate}
\end{defn}

In this definition note that the matrix $\mathbf{X}$ is of rank $s$, then there exists a
rotation matrix $\mathbf{H}\in O(n)$ in order that $\mathbf{HX}=[\mathbf{X}_{1}'|\,\,
\mathbf{0}]'$, where $\mathbf{X}_{1}$ is a non-singular  $s\times s$ matrix, this matrix
absorbs the eulerian meaningful shape of $\mathbf{X}$, and then the matrix $\mathbf{Y}$ is
rotated with that preceding orthogonal matrix, and partitioned according to
$\mathbf{H}\mathbf{Y}=\mathbf{T}=[\mathbf{U}'|\mathbf{V}']'$, which must satisfy
$|(\mathbf{U}'\mathbf{U})(\mathbf{U}'\mathbf{U}+\mathbf{V}'\mathbf{V})^{-1}-\mathbf{I}_{s}|=0$.
Observe also that $\mathbf{U}$ and $\mathbf{V}$ depend on $\mathbf{X}$.

The name of this class of shape obeys the type ${}_{2}F_{1}(\cdot)$ of the hypergeometric
series involved in the main theorem of the present work.

We can put together the above transformations and definitions into the following steps which
obtains the eulerian coordinates $(\mathbf{U}, \mathbf{V})$ of a given pair of original figures
$(\mathbf{X}^{**},\mathbf{Y}^{**})$:
$$\begin{small}
\mathbf{L}\left(%
\begin{array}{cc}
  \mathbf{Y}^{**} & \mathbf{X}^{**} \\
\end{array}%
\right)=\left(%
\begin{array}{cc}
  \mathbf{Y}^{*} & \mathbf{X}^{*} \\
\end{array}%
\right)=\mathbf{W}=\left(%
\begin{array}{cc}
  \mathbf{Y}+\boldsymbol{\mu}_{\mathbf{Y}^{*}} & \mathbf{X}+\boldsymbol{\mu}_{\mathbf{\mathbf{X}}^{*}} \\
\end{array}%
\right)\Rightarrow\left\{%
\begin{array}{ll}
    \mathbf{HX}=\left(%
\begin{array}{c}
  \mathbf{X}_{1} \\
  \textbf{0} \\
\end{array}%
\right)\\
    \mathbf{HY}=\left(%
\begin{array}{c}
  \mathbf{U} \\
  \mathbf{V} \\
\end{array}%
\right)
\end{array}%
\right.,
\end{small}
$$
subject to $
    |(\mathbf{U}'\mathbf{U})(\mathbf{U}'\mathbf{U}+\mathbf{V}'\mathbf{V})^{-1}-\mathbf{I}_{s}|=0$.

Geometrically, Helmertized matrices means that $\mathbf{Y}^{*}$ and $\mathbf{X}^{*}$ are
invariant under translation, and then with the remaining transformations, the resulting
eulerian shape $(\mathbf{U}, \mathbf{V})$ of the pair $(\mathbf{X}^{*},\mathbf{Y}^{*})$ is all
geometrical information remained after filtering out the rotation and permutations of the
columns from the Helmertized pair $(\mathbf{X}^{*},\mathbf{Y}^{*})$.

This last property can be understood as the eulerian shape is invariant under the label of the
axis of coordinates. For example, under this equivalent class, figures in an Euclidian three
dimension space labeled with points with coordinates $(x,y,z)$ are equivalent to figures with
permuted the axis of coordinates, as $(y,z,x)$, $z,x,y$, $(y,x,z)$, $(z,y,x)$ and $(x,z,y)$.
This unusual invariance does not hold in the classical shapes via euclidian, affine or
projective transformations and it is the most important novelty of the proposed approach,
because does not require that the landmarks of the two objects strictly match in labels of the
coordinates of axis.

The few works published in this area belong to procrustes methods based on unlabelled points
where the minimizations of the distances between the figures are carried out over the
similarity transformations and permutations of the landmarks labels, they include a number of
applications as megalithic standing stones ( Cornwall) and human vision, among others, see
section 12.9 of \citet{DM98}. As usual, the problems of inference and the exact density are one
of the obstacles to solve in the existing approaches.  Finally, this invariant statistics of
shape analysis enters too into the so called landmark free approaches, a wider area of study,
see for example section 12.10 of \citet{DM98}.

Clearly, with definition \ref{Def:eulerian} we can implement classical procrustes analysis, as
it was done for Eucliand and affine transformation, but we are interesting in the statistical
analysis of the emergent shape.

In presence of randomness in the location of landmarks in the two Helemertized figures
$\mathbf{X}$ and $\mathbf{Y}$  (linear transformations of the original figures
$\mathbf{X}^{**}$ and $\mathbf{Y}^{**}$), the so called squared canonical correlation
coefficients allow the statistical matching of the objects via
\begin{equation}\label{eulerian2}
    |(\mathbf{Y}'\mathbf{Y})^{-1}\mathbf{Y}'\mathbf{X}(\mathbf{X}'\mathbf{X})^{-1}\mathbf{X}'\mathbf{Y}-r^{2}\mathbf{I}_{K}|=0,
\end{equation}
which is reduced to
\begin{equation}\label{eulerian3}
    |(\mathbf{U}'\mathbf{U})(\mathbf{U}'\mathbf{U}+\mathbf{V}'\mathbf{V})^{-1}-r^2\mathbf{I}_{K}|=0,
\end{equation}
by the items $(1)$ and $(2)$ of Definition \ref{Def:eulerian}.

As we proved, $\mathbf{X}$ and $\mathbf{Y}$ are the $n=N-1\times K$ helmertized matrices, such
that $\mathbf{Y}\sim \mathcal{N}(\mathbf{0},\mathbf{I}_{n}\otimes \mathbf{\Sigma}_{11})$,
$\mathbf{X}\sim \mathcal{N}(\mathbf{0},\mathbf{I}_{n}\otimes \mathbf{\Sigma}_{22})$ and we are
interested in the measure of the relation between $\mathbf{X}$ and $\mathbf{Y}$, where
$\mathbf{\Sigma}_{12}=\mathbf{\Sigma}_{21}'$ is the covariance matrix between the components of
$\mathbf{X}$ and the components of $\mathbf{Y}$. Thus, $\mathbf{Z}=[\mathbf{Y}|\mathbf{X}]\sim
\mathcal{N}(\mathbf{0},\mathbf{I}_{n}\otimes\mathbf{\Sigma})$, where
$$
  \mathbf{\Sigma}=\left(\begin{array}{cc}\mathbf{\Sigma}_{11}&\mathbf{\Sigma}_{12}\\
  \mathbf{\Sigma}_{21}&\mathbf{\Sigma}_{22}\end{array}\right).
$$

The invariant statistics involved in the shape theory context arrives here when we study the
distribution of the so called squared canonical correlation coefficients under the Eulerian
transformation above described, which obtains the squared canonical form of the covariance
matrix as:
$$
\mathbf{\Sigma}=\left(\begin{array}{ccc}
    \mathbf{I}&\mathbf{P}&\mathbf{0}\\
    \mathbf{P}&\mathbf{I}&\mathbf{0}\\
    \mathbf{0}&\mathbf{0}&\mathbf{I}
    \end{array}\right).
$$
where $\mathbf{P}$ is the $K\times K$ diagonal matrix with the squared canonical correlations
$\rho_{1},\ldots,\rho_{K}$ arranged in decreasing order down the diagonal (see \citet{CONST63},
\citet{MR1982}, p. 539).

So, from theorem 11.3.2 of \citet{MR1982}, p. 557 (see also \citet{CONST63}), we have:

Let $\mathbf{A}=\mathbf{Z}'\mathbf{Z}=\left(%
\begin{array}{cc}
  \mathbf{Y}'\mathbf{Y} & \mathbf{Y}'\mathbf{X} \\
  \mathbf{X}'\mathbf{Y} & \mathbf{X}'\mathbf{X} \\
\end{array}%
\right)=\left(%
\begin{array}{cc}
  \mathbf{A}_{11} & \mathbf{A}_{12} \\
  \mathbf{A}_{21} & \mathbf{A}_{22} \\
\end{array}%
\right)$ have the $\mathcal{W}_{2K}(n,\mathbf{\Sigma})$ distribution, where $n\geq 2K$. Then
the joint probability density function of $r_{1}^{2},\ldots,r_{K}^{2}$, the latent roots of
$\mathbf{A}_{11}^{-1}\mathbf{A}_{12}\mathbf{A}_{22}^{-1}\mathbf{A}_{21}$, is
\begin{eqnarray*}
f(r_{1}^{2},\ldots,r_{K}^{2})&=&\prod_{i=1}^{K}(1-\rho_{i}^{2})^{n/2}{}_{2}F_{1}^{(K)}\left(\frac{1}{2}n,
\frac{1}{2}n;\frac{1}{2}K;\mathbf{P}^{2},\mathbf{R}^{2}\right)\\
&&\times
\frac{\pi^{K^{2}/2}}{\Gamma_{K}^{2}\left(\frac{1}{2}K\right)}\frac{\Gamma_{K}\left(\frac{1}{2}n\right)}
{\Gamma_{K}\left[\frac{1}{2}(n-K)\right]}
\\
&&\times
\prod_{i=1}^{K}\left[(r_{i}^2)^{-1/2}(1-r_{i}^{2})^{(n-2K-1)/2}\right]\prod_{i<j}^{K}(r_{i}^{2}-r_{j}^{2}),
\end{eqnarray*}
for $1>r_{1}^{2}>\cdots>r_{1}^{K}>0$, where $\rho_{1}^{2},\ldots,\rho_{K}^{2}$ are the latent
roots of
$\mathbf{\Sigma}_{11}^{-1}\mathbf{\Sigma}_{12}\mathbf{\Sigma}_{22}^{-1}\mathbf{\Sigma}_{21}$,
$\mathbf{P}^{2}=\diag(\rho_{1}^{2},\ldots,\rho_{K}^{2})$ and
$\mathbf{R}^{2}=\diag(r_{1}^{2},\ldots,r_{K}^{2})$.

\section{The main result}\label{sec:polynomialEulerian}

For decades, the above type of expressions could not computed and they forced to asymptotic
studies (see section 11.3.5 of \citet{MR1982}), and it seems that with the recent algorithms of
\citet{KE06} the exact distribution can be used for applications. However, recently
\citet{cdg:14c}, proved that those algoritms do not converge when the inference is performed in
certain shape invariants as the affine case even applied to small size landmark data.

Computation in matrix variate theory via zonal polynomials (see \citet{MR1982}) and invariant
polynomials of several matrix arguments extending the preceding ones (\citet{Davis79},
\citet{D80}) has been the biggest problem in such area, there are hundreds of papers deriving
numbers of theoretical properties involving such polynomials or asymptotic distributions, etc.
but a very few, considering the problem of inference with the exact densities. For decades the
zonal polynomials were impossible to compute until the works of \citet{KE06}, meanwhile the
computation of invariant polynomials of Davis remains unsolved since 1979. However, infinite
series of zonal polynomials has considerable problems and so the performance of inference
associated to the large number of results.

The problem for the computation of one hypergeometric series is addressed by \citet{KE06} in
page 845 : ``Several problems remain open, among them automatic detection of convergence". This
difficulty is amplified when we have to consider the optimization of a large product of such
hypergeometric series, and the addressed algorithms do not allow large truncations to obtain
certain accuracy and the convergence is completely doubtful.

So, in this scenario we have two ways to follow, a computational study which involves also a
machine problem or trying an analytical solution.

Motivated by a similar situation in affine shape studied by \citet{cdg:14c}, this section
attends an analytical solution to the problem of Eulerian shape inference. As before, the
solution resides in the fact that shape densities are usually infinite series of zonal
polynomials, but a few of them can be studied in the setting of the so called polynomial shape
densities, which exploits some characteristics of the parameters of the distributions indexed
by the number of landmarks and dimensions, and by using matrix generalized matrix relations,
the density can be turned into a polynomial. For example, the affine shape theory is
characterized by a confluent hypergeometric function type which can be turned into a polynomial
density by some special domain and using the so called generalized matrix Kummer relation of
\citet{DC15}, a generalization of the well known Kummer relation by \citet{Herz55} with
gaussian kernel. Continuing in this way, we expect that the next type of hypergeometric shape
distribution, the Eulerian of this paper could be transformed into a polynomial by using some
special parametric space and the similar Euler matrix relations. However, the density derived
in the previous section is a modified Eulerian hypergeometric type and involves two matrix
arguments instead the expected one, which facilitates the automatic application of the Euler
matrix relations. Given the connection between the Eulerian shape density and the squared
canonical correlation result of \citet{CONST63}, it seems that  non existence of that relations
explains why the last result was set as an infinite series of products of zonal polynomials,
instead as a polynomial of low degree, as it happens with a number of result based on Gaussian
distributions and Kummer relations.  Then an improvement of the shape density will be applied
in the squared canonical correlation theory too.

Specifically, recall that the hypergeometric function of two $m\times m$ matrices $\mathbf{X}$
and $\mathbf{Y}$ as arguments are given by
\begin{equation}\label{d732}
    {}_{p}F_{q}^{(m)}(a_{1},\ldots,a_{p};b_{1},\ldots,b_{q};\mathbf{X},\mathbf{Y})=\sum_{k=0}^{\infty}\sum_{\kappa}
    \frac{(a_{1})_\kappa\cdots(a_{p})_\kappa}{(b_{1})_\kappa\cdots(b_{q})_\kappa}
    \frac{C_{\kappa}(\mathbf{X})C_{\kappa}(\mathbf{Y})}{k!C_{\kappa}(\mathbf{I}_{m})}.
\end{equation}

The problem is reduced to obtain certain matrix relation which turns certain series of the type
${}_{2}F_{1}^{(m)}(a,b;c;\mathbf{X},\mathbf{Y})$ into an expression of the type
$f(c-a,b;c;\mathbf{X},\mathbf{Y})$ in such way that $c-a$ are negative integers or negative
half-integers, which will turn $f(\cdot)$ into a polynomial of low degree.

The required computable is obtained  from the following matrix relation due to \citet{C14}:

\begin{eqnarray}\label{Euler}
    {}_{2}F_{1}^{(m)}(a_{1},a_{2};b_{1};&&\mathbf{X},\mathbf{Y})=\sum_{k=0}^{\infty}\frac{1}{k!}\sum_{\kappa}
    \frac{(b_{1}-a_{1})_{\kappa}(a_{2})_{\kappa}}{(b_{1})_{\kappa}}
    \\&&\times\int_{O(m)}|\mathbf{I}_{m}-\mathbf{XHYH}'|^{-a_{2}}C_{\kappa}(-\mathbf{XHYH}'
    (\mathbf{I}_{m}-\mathbf{XHYH}')^{-1})(d\mathbf{H})\nonumber;
\end{eqnarray}
if $b_{1}-a_{1}$ is a negative integer, say $-q$, the hypergeometric function is a polynomial
of degree $mq$.

And given that $n=N-1(\geq 2K)$ is the number of helmertized landmarks in $K$ dimension, then
$\frac{K}{2}-\frac{n}{2}=(K-N+1)/2$ is a negative integer or a half negative integer, so we
have arrive to main result of this work:

\begin{thm}\label{th:polynomial}
Let $\mathbf{A}=\mathbf{Z}'\mathbf{Z}=\left(%
\begin{array}{cc}
  \mathbf{Y}'\mathbf{Y} & \mathbf{Y}'\mathbf{X} \\
  \mathbf{X}'\mathbf{Y} & \mathbf{X}'\mathbf{X} \\
\end{array}%
\right)=\left(%
\begin{array}{cc}
  \mathbf{A}_{11} & \mathbf{A}_{12} \\
  \mathbf{A}_{21} & \mathbf{A}_{22} \\
\end{array}%
\right)$ have the $\mathcal{W}_{2K}(n,\mathbf{\Sigma})$ distribution, where $n\geq 2K$. Then
the joint probability density function of $r_{1}^{2},\ldots,r_{K}^{2}$, the latent roots of
$\mathbf{A}_{11}^{-1}\mathbf{A}_{12}\mathbf{A}_{22}^{-1}\mathbf{A}_{21}$, is
\begin{eqnarray*}
f(r_{1}^{2},&&\ldots,r_{K}^{2})=\frac{\pi^{K^{2}/2}}{\Gamma_{K}^{2}\left(\frac{1}{2}K\right)}
\frac{\Gamma_{K}\left(\frac{1}{2}n\right)}{\Gamma_{K}\left[\frac{1}{2}(n-K)\right]}
\prod_{i=1}^{K}(1-\rho_{i}^{2})^{n/2}\prod_{i<j}^{K}(r_{i}^{2}-r_{j}^{2})\\
&&\times \prod_{i=1}^{K}\left[(r_{i}^2)^{-1/2}(1-r_{i}^{2})^{(n-2K-1)/2}\right]
\sum_{l=0}^{\infty}\frac{1}{l!}\sum_{\lambda}\frac{\left(\frac{K}{2}-\frac{n}{2}\right)_{\lambda}
\left(\frac{n}{2}\right)_{\lambda}}{\left(\frac{K}{2}\right)_{\lambda}}
\\&&\times\int_{O(K)}\left|\mathbf{I}_{K}-\mathbf{P}^{2}\mathbf{HR}^{2}\mathbf{H}'\right|^{-n/2}C_{\lambda}
\left(-\mathbf{P}^{2}\mathbf{HR}^{2}\mathbf{H}'(\mathbf{I}_{K}-\mathbf{P}^{2}\mathbf{HR}^{2}\mathbf{H}')^{-1}\right)(d\mathbf{H}),
\end{eqnarray*}
for $1>r_{1}^{2}>\cdots>r_{1}^{K}>0$, where $\rho_{1}^{2},\ldots,\rho_{K}^{2}$ are the latent
roots of
$\mathbf{\Sigma}_{11}^{-1}\mathbf{\Sigma}_{12}\mathbf{\Sigma}_{22}^{-1}\mathbf{\Sigma}_{21}$,
$\mathbf{P}^{2}=\diag(\rho_{1}^{2},\ldots,\rho_{K}^{2})$ and
$\mathbf{R}^{2}=\diag(r_{1}^{2},\ldots,r_{K}^{2})$. If $\frac{K}{2}-\frac{n}{2}$ is a negative
integer, say $-q$, the function is a polynomial of degree $Kq$ in $K$ variables.
\end{thm}

Note that the parametric condition of Eulerian polynomial shape, i.e.
$\frac{K}{2}-\frac{n}{2}$, ($n\geq 2K$) just demands that the objects in  $K$-even(odd)
dimensions are summarized by a number of $N$-odd(even) anatomical landmarks.

Nevertheless the most of the literature of statistical shape analysis is valid for any
dimension $K$, the application of the theory seems to be very restricted to anatomical
landmarks placed on two and three dimensional objects, in fact the classical text on shape
theory, such as \citet{DM98} only uses 2-D landmark data, and 3-D applications are very rarely
given the complexity of the distributions and the registration of invariant critical points.
So, when we consider $K=2$, all the typical landmark data presented for example in \citet{DM98}
can be studied and the integral over the orthogonal group is a feasible and programming task
using the exact formulae for Jack polynomials (which includes real zonal polynomials) due to
\citet{Caro07}.

Next we give some summary of Jack polynomials and exact two dimensional formulae, which
particularized to zonal polynomials, will be needed in some applications.

Let us characterize the Jack symmetric function $J_{\kappa}^{(\alpha)}(y_{1},\ldots,y_{m})$ of
parameter $\alpha$, see \citet{S97}. Consider  a decreasing sequence of nonnegative integers
$\kappa=(k_{1},k_{2},\ldots)$ with only finitely many nonzero terms is said to be a partition
of $k=\sum k_{i}$. Let $\kappa$ and $\lambda=(l_{1},l_{2},\ldots)$ be two partitions of $k$. We
write $\lambda\leq\kappa$ if $\sum_{i=1}^{t}l_{i}\leq\sum_{i=1}^{t}k_{i}$ for each $t$. The
conjugate of $\kappa$ is $\kappa'=(k_{1}',k_{2}',\ldots)$ where $k_{i}'=
\mbox{card}\{j:k_{j}\geq i\}$. The length of $\kappa$ is $l(k)=\mbox{max}\{i:k_{i}\neq
0\}=k_{1}'$.  If $l(\kappa)\leq m$, one often writes $\kappa=(k_{1},k_{2},\ldots,k_{m})$.  The
partition $(1,\ldots,1)$ of length $m$ will be denoted by $1_{m}$.

And recall that the monomial symmetric function $M_{\kappa}(\cdot)$ indexed by a partition
$\kappa$ can be regarded as a function of an arbitrary number of variables such that all but a
finite number are equal to $0$:  if $y_{i}=0$ for $i>m\geq l(\kappa)$ then
$M_{\kappa}(y_{1},\ldots,y_{m})=\sum y_{1}^{\sigma_{1}}\cdots y_{m}^{\sigma_{m}}$, where the
sum is over all distinct permutations $\{\sigma_{1},\ldots,\sigma_{m}\}$ of
$\{k_{1},\ldots,k_{m}\}$, and if $l(\kappa)>m$ then $M_{\kappa}(y_{1},\ldots,y_{m})=0$.  A
symmetric function $f$ is a linear combination of monomial symmetric functions.  If $f$ is a
symmetric function then $f(y_{1},\ldots,y_{m},0)=f(y_{1},\ldots,y_{m})$.  For each $m\geq 1$,
$f(y_{1},\ldots,y_{m})$ is a symmetric polynomial in $m$ variables.

Thus the Jack symmetric function
$J_{\kappa}^{(\alpha)}(y_{1},\ldots,y_{m})$ with a parameter
$\alpha$, satisfy the following conditions:
\begin{eqnarray}
    J_{\kappa}^{(\alpha)}(y_{1},\ldots,y_{m})=\sum_{\lambda\leq\kappa}j_{\kappa,\lambda}
    M_{\lambda}(y_{1},\ldots,y_{m}),\label{eq:Condition1}\\
    J_{\kappa}^{(\alpha)}(1,\ldots,1)=\alpha^{k}\prod_{i=1}^{m}\left(\frac{m-i+1}{\alpha}\right)_{k_{i}},
    \label{eq:Condition2}\\
    \sum_{i =1}^{m}y_{i}^{2}\frac{\partial^{2}J_{\kappa}^{(\alpha)}(y_{1},\ldots,y_{m})}{\partial y_{i}^{2}}
    +\frac{2}{\alpha}\sum_{i=1}^{m}y_{i}^{2}\sum_{j\neq i}\frac{1}{y_{i}-y_{j}}
    \frac{\partial J_{\kappa}^{(\alpha)}(y_{1},\ldots,y_{m})}{\partial y_{i}}=\nonumber\\
    \sum_{i=1}^{m}k_{i}(k_{i}-1+\frac{2}{\alpha}(m-i))J_{\kappa}^{(\alpha)}(y_{1},\ldots,y_{m})
    \label{eq:Condition3}.
\end{eqnarray}
Here the constants $j_{\kappa,\lambda}$ do not dependent on
$y_{i}'s$ but on $\kappa$ and $\lambda$, and
$(a)_{n}=\prod_{i=1}^{n}(a+i-1)$. Note that if $m<l(\kappa)$ then
$J_{\kappa}^{(\alpha)}(y_{1},\ldots,y_{m})=0$. The conditions
include the case $\alpha =0$ and then
$J_{\kappa}^{(0)}(y_{1},\ldots,y_{m})=e_{\kappa
'}\prod_{i=1}^{m}(m-i+1)^{k_{i}}$, where
$e_{\kappa}(y_{1},\ldots,y_{m})=\prod_{i=1}^{l(\kappa)}e_{k_{i}}(y_{1},\ldots,y_{m})$
are the elementary symmetric functions indexed by partitions
$\kappa$, if $m\geq l(\kappa)$ then
$e_{r}(y_{1},\ldots,y_{m})=\sum_{i_{1}<i_{2}<\cdots<i_{r}}y_{i_{1}}\cdots
y_{i_{r}}$, and if $m< l(\kappa)$ then
$e_{r}(y_{1},\ldots,y_{m})=0$, see \citet{S97}.

Now, from \citet{KE06}, the Jack functions $J_{\kappa}^{(\alpha)}(\mathbf{Y}) =
J_{\kappa}^{(\alpha)}(y_{1},\ldots,y_{m})$, with $y_{1},\ldots,y_{m}$ being the eigenvalues of
the matrix $\mathbf{Y}$, can be normalized in such a way that
$$
  \sum_{\kappa}C_{\kappa}^{\alpha}(\mathbf{Y}) = (\tr(\mathbf{Y}))^{k},
$$
where $C_{\kappa}^{\alpha}(\mathbf{Y})$ denotes the Jack polynomials. They are related to the
Jack functions by
\begin{equation}\label{eq1}
    C_{\kappa}^{\alpha}(\mathbf{Y})= \frac{\alpha^{k} k!}{j_{\kappa}} J_{\kappa}^{\alpha}(\mathbf{Y}),
\end{equation}
where
$$
  j_{\kappa} = \prod_{(i,j) \in \kappa} h_{*}^{\kappa}(i,j) h^{*}_{\kappa}(i,j),
$$
and $h_{*}^{\kappa}(i,j) = k_{j}-i+\alpha(k_{i}-j+1)$ and
$h^{*}_{\kappa}(i,j) = k_{j}-i+1+\alpha(k_{i}-j)$ are the upper and
lower hook lengths at $(i,j) \in \kappa$, respectively.

Then by applying (\ref{eq1}), we can write (\ref{eq:Condition3}) as
$$
  \sum_{1}^{m}y_{i}^{2}\frac{\partial^{2}C_{\kappa}^{(\alpha)}(\mathbf{Y})}{\partial y_{i}^{2}}
    +\frac{2}{\alpha}\sum_{i=1}^{m}y_{i}^{2}\sum_{j\neq i}\frac{1}{y_{i}-y_{j}}
    \frac{\partial C_{\kappa}^{(\alpha)}(\mathbf{Y})}{\partial y_{i}} =  \hspace{5cm}
$$
\begin{equation}\label{eqdif}
    \hspace{6cm}\sum_{i=1}^{m}k_{i}(k_{i}-1+\frac{2}{\alpha}(m-i))C_{\kappa}^{(\alpha)}(\mathbf{Y}).
\end{equation}
Now, when $m=2$ in (\ref{eqdif}), \citet{Caro07} found the following
formulae for Jack Polynomials of the Second Order
\small{%
$$\hspace{-1cm}
  \displaystyle\frac{C^{(\alpha)}_{(k_{1},k_{2})}(\mathbf{Y})}{C^{(\alpha)}_{(k_{1},k_{2})}(\mathbf{I}_{2})} =
    (y_{1}y_{2})^{(k_{1}+k_{2})/2}A_{1}F\left(-\frac{\rho}{2},\frac{\rho}{2}+\frac{1}{\alpha};\frac{1}{2};
    \frac{(y_{1}+y_{2})^2}{4y_{1}y_{2}}\right)
$$}
\begin{equation}\label{solgeneral}
    + \frac{(y_{1}y_{2})^{(k_{1}+k_{2}-1)/2}}{2(y_{1}+y_{2})^{-1}}
    A_{2}
     F\left(\frac{1}{\alpha} + \frac{1+\rho}{2},\frac{1}{2} -
   \frac{\rho}{2};\frac{3}{2};\frac{(y_{1}+y_{2})^2}{4y_{1}y_{2}}\right),
\end{equation}
with $\rho$ being either even or odd. For distinguishing the case
under consideration, odd or even, we will use the upper indices $o$
or $e$ with $A_{1}$ and $A_{2}$. Then the corresponding solutions
are the following

\textbf{Even case.} If $\rho=k_{1}-k_{2}=2n$, $n=0,1,2,\ldots$ then
$$
  A_{1}^{e}= \frac{(-1)^{n} \displaystyle\prod_{i=0}^{n-1}(1+2i)}
    {\displaystyle\prod_{i=0}^{n-1}\left(1+2\left(\frac{1}{\alpha}+i\right)\right)}
    \quad \mbox{and}\quad A_{2}^{e}= 0.
$$

\textbf{Odd case.} If $\rho=k_{1}-k_{2}=2n +1$, $n=0,1,2,\ldots$
then
$$
  A_{1}^{o} = 0 \quad \mbox{and}\quad A_{2}^{o}= (2n +1) A_{1}^{e}.
$$

Three particular cases are of interest in the literature: the
quaternionic case ($\alpha = 1/2$), the complex zonal polynomials
($\alpha =1 $) and the real zonal polynomials ($\alpha = 2$), these
results are summarized in the following table:

\bigskip
\begin{tabular}{||c||c|c|c|c|c|c||}
  \hline\hline
  &&&&&&\\
  $\alpha$ & $\rho$ & $a$ & $b$ & $c$ & $A_{1}$ & $A_{2}$\\
  \hline\hline
  \multicolumn{1}{||c||}{}&&&&&& \\
  $\displaystyle\frac{1}{2}$ & even & $-n$ & $n +2$ & $\displaystyle\frac{1}{2}$ &
    $\displaystyle\frac{(-1)^{n}3}{(2n + 1)(2n+3)}$ &$0$\\
  \multicolumn{1}{||c||}{}&&&&&& \\
  \cline{2-7}
  \multicolumn{1}{||c||}{}&&&&&& \\
   & odd & $n + 3$ & $-n$ & $\displaystyle\frac{3}{2}$ & 0 & $\displaystyle\frac{(-1)^{n}3}{(2n+3)}$\\
 \multicolumn{1}{||c||}{}&&&&&& \\
 \hline\hline
 \multicolumn{1}{||c||}{}&&&&&& \\
  $1$ & even & $-n$ & $n +1$ & $\displaystyle\frac{1}{2}$ & $\displaystyle\frac{(-1)^{n}}{(2n + 1)}$ &$0$\\
  \multicolumn{1}{||c||}{}&&&&&& \\
  \cline{2-7}
  \multicolumn{1}{||c||}{}&&&&&& \\
   & odd & $n + 2$ & $-n$ & $\displaystyle\frac{3}{2}$ & 0 & $(-1)^{n}$\\
 \multicolumn{1}{||c||}{}&&&&&& \\
 \hline\hline
 \multicolumn{1}{||c||}{}&&&&&& \\
  $2$ & even & $-n$ & $n +1/2$ & $\displaystyle\frac{1}{2}$ &
    $\displaystyle\frac{(-1)^{n}(2n)!}{2^{2n}(n!)^{2}}$ &$0$\\
  \multicolumn{1}{||c||}{}&&&&&& \\
  \cline{2-7}
  \multicolumn{1}{||c||}{}&&&&&& \\
   & odd & $n + 3/2$ & $-n$ & $\displaystyle\frac{3}{2}$ & 0 & $\displaystyle\frac{(-1)^{n}(2n+1)!}{2^{2n}(n!)^{2}}$\\
 \multicolumn{1}{||c||}{}&&&&&& \\
 \hline\hline
\end{tabular}

\bigskip

The above formula for the real zonal polynomials corresponds to that derived by \citet{JAT68}
and for the complex zonal polynomials, obtained by \citet{Caro06b}, meanwhile it is the first
appearance of an exact formulae for quaternionic polynomials.

So, finally, for the integration demanded in Theorem \ref{th:polynomial} over the orthogonal
group $O(2)$ with respect the normalized invariant measure $(d\mathbf{H})$, we just follow the
usual parametrization techniques, see for example \citet{MR1982}.

\section{Applications}\label{sec:applications}
In this section we describe the general procedure for applications of the Eulerian shape
distributions and compare the computational performance when the inference of a certain
appropriate landmark data for this data is studied under the modified classical result of
\citet{CONST63} via the algorithms of \citet{KE06} and the exact polynomial distribution here
derived.

For the first aspect, there are a number of situations to consider, for example, the Eulerian
shape can be used to establish correlation structure of two populations, or determine the
correlation structure of a population and a given template, also, the correlation structure of
a population and a given object probably belonging to another population, and finally, we can
use the polynomial distributions for exploring the problem of landmark discrimination. Next we
explain the associated procedure of each one.

\subsection{Correlation structure of two populations}

Consider a sample of $m$ figures of a population $\mathbf{X}$ and $m$ figures of a population
$\mathbf{Y}$; both samples satisfying the shape invariance assumed for the Eulerian shape and
the fact that $(K-N+1)/2$, ($N-1\geq 2K$) is a negative integer, i.e. the objects in
$K$-even(odd) dimensions are summarized by a number of $N$-odd(even) anatomical landmarks. For
each pair of figures $\mathbf{Z}=[\mathbf{Y}|\mathbf{X}]$ compute the corresponding
$r_{1,j}^{2},\ldots,r_{K,j}^{2}$, $j=1,\ldots,m$ the latent roots of
$\mathbf{A}_{11}^{-1}\mathbf{A}_{12}\mathbf{A}_{22}^{-1}\mathbf{A}_{21}$, in order to obtain
the function $f_{j}(r_{1}^{2},\ldots,r_{K}^{2})$ of Theorem \ref{th:polynomial}.

Then maximize the likelihood function
$$
L(\rho_{1}^{2},\ldots,\rho_{K}^{2})=\prod_{j=1}^{m}f_{j}(r_{1}^{2},\ldots,r_{K}^{2})
$$
respect to the population squared canonical correlation coefficients
$\rho_{1}^{2},\ldots,\rho_{K}^{2}$.

These estimates will provide the strength of relationship between the two populations in the
sense of the Eulerian shape invariance.

As we shall see in the experiments various methods can be implemented in order to obtain more
information about the shape equality of both populations without using the classical statistics
for means, instead we can use the exact distribution of certain template vs one of the
population and obtain the probability than an extreme correlation structure of the other
population influenced by the same template, can be occurred.

\subsection{Correlation structure of a population and a given
template}

Exactly as before, but in this case, we have $m$ repeated templates $\mathbf{X}$ and we want to
verify the  correlation between the population $\mathbf{Y}$ and the template $\mathbf{X}$ under
the Eulerian shape assumptions.

\subsection{Correlation structure of a population and a given
figure}

Exactly as before, but in this case, we have $m$ repeated figures $\mathbf{X}$ and we want to
verify the correlation between the population $\mathbf{Y}$ and a given ``figure" $\mathbf{X}$,
i.e. it can explain, under the requirements of the Eulerian invariants,  if the last one
belongs  or not to the population $\mathbf{Y}$.

\subsection{Landmark discrimination}

This is one of the main less studied problems in the shape theory. It concerns the reduction of
set of landmarks for a given figure without changing the explanation of the object. In other
words it worries about reduction of dimensionality of a figure and set the shape theory in the
domain of very low dimensions. The method of the preceding section can be used in order to
study the landmark discrimination. To get this end, we start with figures $\mathbf{Y}$ and
$\mathbf{X}$ (in any of the three possible scenarios explained in the above items) with a large
number of landmarks $N$ (not necessary anatomical landmarks, they can be mathematical or
pseudolandmarks), estimate the corresponding $\rho's$, then reduce the number of landmarks and
repeat the estimation until the $\rho's$ change drastically, then we reach an optimal number of
landmarks which summaries properly the object.

This technique is particular appropriate for shape analysis involving mathematical or
pseudomathematical landmarks, which usually are equally-placed in the contour of the object.
For example, the method is applicable in some extreme cases, such as discrimination of classes
of potatoes or differentiation of non regular objects with no chance of placing anatomical
landmarks objets. A profuse study of this landmark discrimination methodology is will be
carried out in a subsequent work. Another useful application under consideration tries an
automatic detection of anatomical landmarks using their limit from pseudo critical points. For
example, the 46 pseudo-landmarks defined by \citet{DM98} in the contour of the mouse vertebrae
can be used for studying the discrimination problem and finding the possible optimal set of six
anatomical landmarks usually referred in literature.

\subsection{An application in postcode recognition}

We end this section by illustrating the addressed problem of inference with the distribution
provided by \citet{CONST63} and the algorithms for hypergeometric functions by \citet{KE06}.
Then we compare the above results with the polynomial density here derived.

We select a known landmark data of literature on postcode recognition which has presented
problems in inference using another shape invariants as the affine shape, see \citet{cdg:14c}.

\begin{figure}[!h]
  \begin{center}
  \includegraphics[width=13cm,height=8cm]{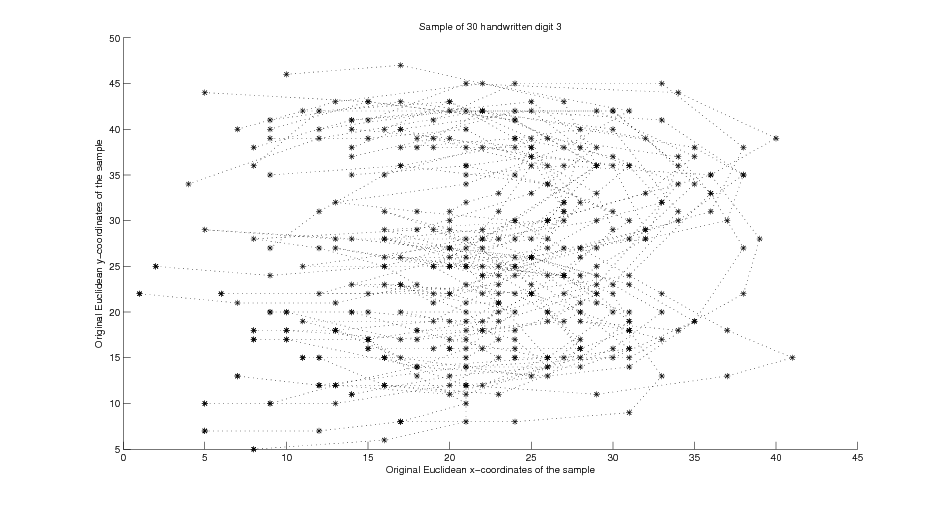}
  \caption{30 handwritten digit 3.}\label{fig:sample}
  \end{center}
\end{figure}

\begin{figure}[!h]
  \begin{center}
  \includegraphics[width=13cm,height=8cm]{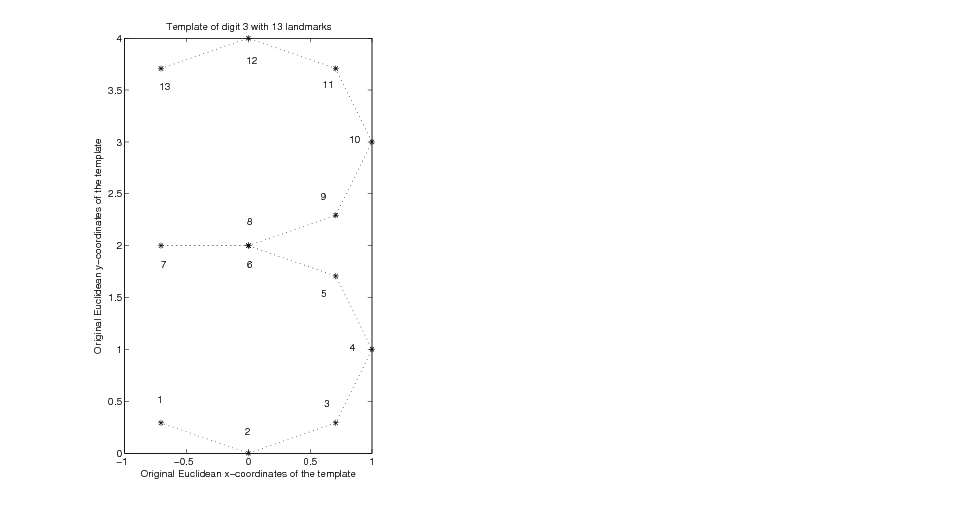}
  \caption{Configuration template.}\label{fig:template}
  \end{center}
\end{figure}

In this case, we have a 30 sample of handwritten digit 3 (see Figure \ref{fig:sample}), with
$N=13, (K=2)$ landmarks placed according an equally spaced two joint octagons as a template
(see Figure \ref{fig:template}).  We are interested in estimating the strength the estimate of
correlation structure between the template and the handwritten digit 3. The main theorem can be
used to obtain the estimates of population squared canonical correlations $\rho_{1}^{2}$ and
$\rho_{2}^{2}$ by optimizing the corresponding likelihood.  In this case the distribution is a
polynomial of degree $-\frac{1}{2}(K-N+1)\times K=10$ in the two latent roots of the matrix
argument in the zonal polynomial.

If we truncate the original distributions derived by \citet{CONST63} and use the algorithms of
\citet{KE06} we obtain the following table, which shows the starting point
$(r_{1i}^{2},r_{2i}^{2})$for the optimization, the population estimate coefficients between the
template and the handwritten digit 3, the convergence and the number of iteractions of the
algorithm. In this case we just truncate the infinite series of zonal polynomials at 10th
degree in order to compare the performance with the polynomial of 10-th degree which provides
the solution. The computations were performed with a processor Intel(R) Corel(TM)2 Duo CPU,
E7400@2.80GHz, and 2,96GB of RAM.

\begin{scriptsize}
\begin{tabular}{||c|c|c|c|c|c|c||}
  \hline
  $r_{1i}^{2}$ &$r_{2i}^{2}$&Fval.&Conv.&Iter.&$\widehat{\rho}_{1}^{\,\,2}$&$\widehat{\rho}_{2}^{\,\,2}$  \\\hline\hline
     0.1 & 0.1 & 372,8472 & y & 23 & 0.1260 & 0.1260 \\
0.1 & 0.2 & 372.8472 & y & 23 & 0.1260 & 0.1260 \\
0.1 & 0.3 & 372.8472 & y & 27 & 0.1260 & 0.1260 \\
0.1 & 0.4 & 372.8472 & y & 30 & 0.1260 & 0.1259 \\
0.1 & 0.5 & 372.8472 & y & 44 & 0.1260 & 0.1260 \\
0.1 & 0.6 & 372.8472 & y & 31 & 0.1260 & 0.1260 \\
0.1 & 0.7 & 372.8472 & y & 31 & 0.1260 & 0.1260 \\
0.1 & 0.8 & 372.8472 & y & 32 & 0.1260 & 0.1260 \\
0.1 & 0.9 & - & n & - & - & - \\
 \hline
\end{tabular}
\end{scriptsize}
\begin{scriptsize}
\begin{tabular}{||c|c|c|c|c|c|c||}
  \hline
  $r_{1i}^{2}$ &$r_{2i}^{2}$&Fval.&Conv.&Iter.&$\widehat{\rho}_{1}^{\,\,2}$&$\widehat{\rho}_{2}^{\,\,2}$  \\\hline\hline
0.2 & 0.2 & 372.8472 & y & 27 & 0.1260 & 0.1260 \\
0.2 & 0.3 & 372.8472 & y & 34 & 0.1259 & 0.1260 \\
0.2 & 0.4 & 372.8472 & y & 32 & 0.1259 & 0.1259 \\
0.2 & 0.5 & 372.8472 & y & 37 & 0.1259 & 0.1260 \\
0.2 & 0.6 & 372.8472 & y & 36 & 0.1260 & 0.1260 \\
0.2 & 0.7 & - & n & - & - & - \\
0.2 & 0.8 & - & n & - & - & - \\
0.2 & 0.9 & - & n & - & - & - \\
 \hline
\end{tabular}
\end{scriptsize}
\begin{scriptsize}
\begin{tabular}{||c|c|c|c|c|c|c||}
  \hline
  $r_{1i}^{2}$ &$r_{2i}^{2}$&Fval.&Conv.&Iter.&$\widehat{\rho}_{1}^{\,\,2}$&$\widehat{\rho}_{2}^{\,\,2}$  \\\hline\hline
0.3 & 0.3 & 372.8472 & y & 29 & 0.1260 & 0.1260 \\
0.3 & 0.4 & 372.8472 & y & 39 & 0.1260 & 0.1260 \\
0.3 & 0.5 & 372.8472 & y & 40 & 0.1260 & 0.1260 \\
0.3 & 0.6 & - & n & - & - & - \\
0.3 & 0.7 & - & n & - & - & - \\
0.3 & 0.8 & - & n & - & - & - \\
0.3 & 0.9 & - & n & - & - & - \\
 \hline
\end{tabular}
\end{scriptsize}
\begin{scriptsize}
\begin{tabular}{||c|c|c|c|c|c|c||}
  \hline
  $r_{1i}^{2}$ &$r_{2i}^{2}$&Fval.&Conv.&Iter.&$\widehat{\rho}_{1}^{\,\,2}$&$\widehat{\rho}_{2}^{\,\,2}$  \\\hline\hline
0.4 & 0.4 & - & n & - & - & - \\
0.4 & 0.5 & 372.8472 & y & 40 & 0.1260 & 0.1260 \\
0.4 & 0.6 & - & n & - & - & - \\
0.4 & 0.7 & - & n & - & - & - \\
0.4 & 0.8 & - & n & - & - & - \\
0.4 & 0.9 & - & n & - & - & - \\
 \hline
\end{tabular}
\end{scriptsize}
\begin{scriptsize}
\begin{tabular}{||c|c|c|c|c|c|c||}
  \hline
  $r_{1i}^{2}$ &$r_{2i}^{2}$&Fval.&Conv.&Iter.&$\widehat{\rho}_{1}^{\,\,2}$&$\widehat{\rho}_{2}^{\,\,2}$  \\\hline\hline
0.5 & 0.5 & - & n & - & - & - \\
0.5 & 0.6 & 372.8472 & y & 37 & 0.1260 & 0.1260 \\
0.5 & 0.7 & - & n & - & - & - \\
0.5 & 0.8 & - & n & - & - & - \\
0.5 & 0.9 & - & n & - & - & - \\
 \hline
\end{tabular}
\end{scriptsize}
\begin{scriptsize}
\begin{tabular}{||c|c|c|c|c|c|c||}
  \hline
  $r_{1i}^{2}$ &$r_{2i}^{2}$&Fval.&Conv.&Iter.&$\widehat{\rho}_{1}^{\,\,2}$&$\widehat{\rho}_{2}^{\,\,2}$  \\\hline\hline
0.6 & 0.6 & 372.8472 & y & 34 & 0.1259 & 0.1260 \\
0.6 & 0.7 & 372.8472 & y & 37 & 0.1260 & 0.1260 \\
0.6 & 0.8 & - & n & - & - & - \\
0.6 & 0.9 & - & n & - & - & - \\
 \hline
\end{tabular}
\end{scriptsize}
\begin{scriptsize}
\begin{tabular}{||c|c|c|c|c|c|c||}
  \hline
  $r_{1i}^{2}$ &$r_{2i}^{2}$&Fval.&Conv.&Iter.&$\widehat{\rho}_{1}^{\,\,2}$&$\widehat{\rho}_{2}^{\,\,2}$  \\\hline\hline
0.7 & 0.7 & - & n & - & - & - \\
0.7 & 0.8 & - & n & - & - & - \\
0.7 & 0.9 & - & n & - & - & - \\
0.8 & 0.8 & - & n & - & - & - \\
0.8 & 0.9 & - & n & - & - & - \\
0.9 & 0.9 & - & n & - & - & - \\
 \hline
\end{tabular}
\end{scriptsize}

There is no doubt that the algorithm does not work and cannot reach any credible solution; in
fact for some starting points the algorithm does not converge even at a very low truncation of
10th degree.

Meanwhile the polynomial density gives a fast convergence no matters the starting point we
consider, the result in this case by different classical methods for optimization, is:
$\widehat{\rho}_{1}^{\,\,2}=0.9542$ and $\widehat{\rho}_{2}^{\,\,2}=0.6740$.

Now, if we want to approximate the exact result by using the result of Muihead, the  series
must be truncated in a large number as the following table shows. We highlight that the
algorithm was intervened many times in the starting point in order to obtain an stable
solution.

\begin{scriptsize}
\begin{tabular}{||c|c|c||||c|c|c||||c|c|c||}
  \hline
  Truncation &$\widehat{\rho}_{1}^{\,\,2}$&$\widehat{\rho}_{2}^{\,\,2}$&  Truncation &$\widehat{\rho}_{1}^{\,\,2}$
  &$\widehat{\rho}_{2}^{\,\,2}$ &Truncation &$\widehat{\rho}_{1}^{\,\,2}$&$\widehat{\rho}_{2}^{\,\,2}$\\\hline\hline
  20& 0.5890&    0.5889  &120&    0.9246&    0.6697     &220  &0.9484&    0.6733\\
  40& 0.7319 &   0.7318&140   &0.9330&    0.6710&240  &0.9500 &   0.6735\\
  60  &0.8554 &   0.6713&160  &0.9389 &   0.6719&260  &0.9512  &  0.6736\\
  80  &0.8922  &  0.6670&180  &0.9431  &  0.6725&280  &0.9521   & 0.6736\\
  100 &0.9121   & 0.6682&200& 0.9461&    0.6730 & 300    &0.9527&    0.6737\\
 \hline\hline
\end{tabular}
\end{scriptsize}

 \begin{scriptsize}
\begin{tabular}{||c|c|c||||c|c|c||}
  \hline
  Truncation &$\widehat{\rho}_{1}^{\,\,2}$&$\widehat{\rho}_{2}^{\,\,2}$&Truncation &$\widehat{\rho}_{1}^{\,\,2}$&$\widehat{\rho}_{2}^{\,\,2}$
    \\\hline\hline
320 &0.9531 &   0.6737&420  &0.9538  &  0.6737\\
340  &0.9534 &   0.6737&440 &0.9538   & 0.6737\\
360 &0.9536   & 0.6737&460  &0.9539&    0.6736\\
380 &0.9537&    0.6737&480  &0.9539 &   0.6737\\
400 &0.9538 &   0.6736&500  &0.9539   & 0.6737\\
 \hline\hline
\end{tabular}
\end{scriptsize}

With the implemented routines based on infinite series of zonal polynomials was impossible to
reach the solution obtained with the exact distribution. After extremely large truncation of
500 they could not stabilize to the solution obtained with exact polynomial of 10-th degree.
Note, that we cannot provide an analytical of numerical relationship  between the truncation
and the estimates, it is exactly the open problem addressed by \citet{KE06}.

Once the computation problem is solved, we must pay attention to the interpretation of the
squared canonical correlation estimates, $\widehat{\rho}_{1}^{\,\,2}=0.9542$ and
$\widehat{\rho}_{2}^{\,\,2}=0.6740$, corresponding to the comparison between the template and a
typical handwritten sample of digit 3. The first approach consist of computing the probability
that a given handwritten digit 3 with particular $r_{t1}^{2}$ and $r_{t2}^{2}$ is far from the
template, in this case we just need to compute
$P((r_{1}^{2},r_{2}^{2})>(r_{t1}^{2},r_{t2}^{2}))$, which is a feasible task by using the
polynomial densities of Theorem \ref{th:polynomial}.

\begin{figure}[!h]
  \begin{center}
  \includegraphics[width=13cm,height=8cm]{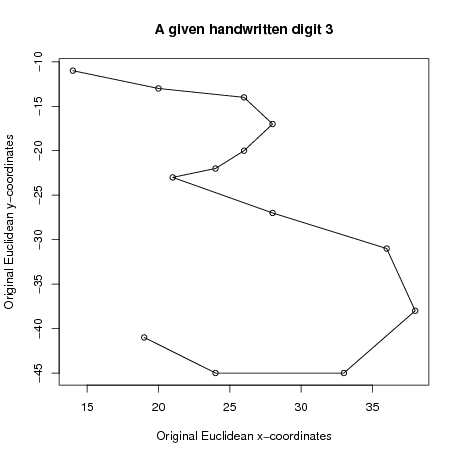}
  \caption{A randomly selected handwritten digit 3.}\label{fig:digit3sample}
  \end{center}
\end{figure}

Now, to obtain a measure for interpretation, we can apply the same procedure to a non-template
object, for example, any handwritten digit 3 selected at random, see Figure
\ref{fig:digit3sample}; the next table summarizes the convergent problem of optimization when
we use the infinite series. First we demonstrate the inefficient of the series which is
truncated at the level of the degree 10-th of the exact polynomial distribution.

\begin{scriptsize}
\begin{tabular}{||c|c|c|c|c|c|c||}
  \hline
  $r_{1i}^{2}$ &$r_{2i}^{2}$&Fval.&Conv.&Iter.&$\widehat{\rho}_{1}^{\,\,2}$&$\widehat{\rho}_{2}^{\,\,2}$  \\\hline\hline
     0.1 & 0.1 & 553.8395
 & y & 24 & 0.1283 & 0.1282 \\
0.1 & 0.2 & 553.8395 & y & 24 & 0.1283 & 0.1282 \\
0.1 & 0.3 & 553.8395 & y & 30 & 0.1282 & 0.1282 \\
0.1 & 0.4 & 553.8395 & y & 33 & 0.1282 & 0.1283 \\
0.1 & 0.5 & 553.8395 & y & 35 & 0.1282 & 0.1282 \\
0.1 & 0.6 & 553.8395 & y & 32 & 0.1282 & 0.1282 \\
0.1 & 0.7 & 553.8395 & y & 31 & 0.1282 & 0.1282 \\
0.1 & 0.8 & 553.8395 & y & 34 & 0.1283 & 0.1282 \\
0.1 & 0.9 & - & n & - & - & - \\
 \hline
\end{tabular}
\end{scriptsize}
\begin{scriptsize}
\begin{tabular}{||c|c|c|c|c|c|c||}
  \hline
  $r_{1i}^{2}$ &$r_{2i}^{2}$&Fval.&Conv.&Iter.&$\widehat{\rho}_{1}^{\,\,2}$&$\widehat{\rho}_{2}^{\,\,2}$  \\\hline\hline
0.2 & 0.2 & 553.8395 & y & 26 & 0.1282 & 0.1282 \\
0.2 & 0.3 & 553.8395 & y & 34 & 0.1282 & 0.1282 \\
0.2 & 0.4 & 553.8395 & y & 33 & 0.1282 & 0.1282 \\
0.2 & 0.5 & 553.8395 & y & 38 & 0.1282 & 0.1282 \\
0.2 & 0.6 & 553.8395 & y & 38 & 0.1282 & 0.1283 \\
0.2 & 0.7 & - & n & - & - & - \\
0.2 & 0.8 & - & n & - & - & - \\
0.2 & 0.9 & - & n & - & - & - \\
 \hline
\end{tabular}
\end{scriptsize}
\begin{scriptsize}
\begin{tabular}{||c|c|c|c|c|c|c||}
  \hline
  $r_{1i}^{2}$ &$r_{2i}^{2}$&Fval.&Conv.&Iter.&$\widehat{\rho}_{1}^{\,\,2}$&$\widehat{\rho}_{2}^{\,\,2}$  \\\hline\hline
0.3 & 0.3 & 553.8395 & y & 30 & 0.1282 & 0.1282 \\
0.3 & 0.4 & 553.8395 & y & 37 & 0.1282 & 0.1283 \\
0.3 & 0.5 & 553.8395 & y & 37 & 0.1282 & 0.1283 \\
0.3 & 0.6 & 553.8395 & y & 36 & 0.1283 & 0.1282 \\
0.3 & 0.7 & - & n & - & - & - \\
0.3 & 0.8 & - & n & - & - & - \\
0.3 & 0.9 & - & n & - & - & - \\
 \hline
\end{tabular}
\end{scriptsize}
\begin{scriptsize}
\begin{tabular}{||c|c|c|c|c|c|c||}
  \hline
  $r_{1i}^{2}$ &$r_{2i}^{2}$&Fval.&Conv.&Iter.&$\widehat{\rho}_{1}^{\,\,2}$&$\widehat{\rho}_{2}^{\,\,2}$  \\\hline\hline
0.4 & 0.4 & - & n & - & - & - \\
0.4 & 0.5 & 553.8395 & y & 39 & 0.1282 & 0.1282 \\
0.4 & 0.6 & - & n & - & - & - \\
0.4 & 0.7 & - & n & - & - & - \\
0.4 & 0.8 & - & n & - & - & - \\
0.4 & 0.9 & - & n & - & - & - \\
 \hline
\end{tabular}
\end{scriptsize}
\begin{scriptsize}
\begin{tabular}{||c|c|c|c|c|c|c||}
  \hline
  $r_{1i}^{2}$ &$r_{2i}^{2}$&Fval.&Conv.&Iter.&$\widehat{\rho}_{1}^{\,\,2}$&$\widehat{\rho}_{2}^{\,\,2}$  \\\hline\hline
0.5 & 0.5 & - & n & - & - & - \\
0.5 & 0.6 & 553.8395 & y & 36 & 0.1282 & 0.1282 \\
0.5 & 0.7 & - & n & - & - & - \\
0.5 & 0.8 & - & n & - & - & - \\
0.5 & 0.9 & - & n & - & - & - \\
 \hline
\end{tabular}
\end{scriptsize}
\begin{scriptsize}
\begin{tabular}{||c|c|c|c|c|c|c||}
  \hline
  $r_{1i}^{2}$ &$r_{2i}^{2}$&Fval.&Conv.&Iter.&$\widehat{\rho}_{1}^{\,\,2}$&$\widehat{\rho}_{2}^{\,\,2}$  \\\hline\hline
0.6 & 0.6 & - & n & - & - & - \\
0.6 & 0.7 & 553.8395 & y & 40 & 0.1282 & 0.1282 \\
0.6 & 0.8 & - & n & - & - & - \\
0.6 & 0.9 & - & n & - & - & - \\
 \hline
\end{tabular}
\end{scriptsize}
\begin{scriptsize}
\begin{tabular}{||c|c|c|c|c|c|c||}
  \hline
  $r_{1i}^{2}$ &$r_{2i}^{2}$&Fval.&Conv.&Iter.&$\widehat{\rho}_{1}^{\,\,2}$&$\widehat{\rho}_{2}^{\,\,2}$  \\\hline\hline
0.7 & 0.7 & - & n & - & - & - \\
0.7 & 0.8 & - & n & - & - & - \\
0.7 & 0.9 & - & n & - & - & - \\
0.8 & 0.8 & - & n & - & - & - \\
0.8 & 0.9 & - & n & - & - & - \\
0.9 & 0.9 & - & n & - & - & - \\
 \hline
\end{tabular}
\end{scriptsize}

Again the malfunction of the algorithms is revealed in the simple truncation 10; no convergence
is obtained for every starting point, in particular for the largest ones which are near to the
mean. Each pair of solutions is so far from the estimates based on the exact density; those
values were derived by different numerical methods, no problems of convergence was occurred and
were independent of the starting points, the corresponding estimates are:
$\widehat{\rho}_{1}^{\,\,2}=0.9753$ and $\widehat{\rho}_{2}^{\,\,2}=0.8748$.

If we want to pursuit that solution we require enormous values for truncating the series, as it
is shown in the following table.

\begin{scriptsize}
\begin{tabular}{||c|c|c||||c|c|c||||c|c|c||}
  \hline
  Truncation &$\widehat{\rho}_{1}^{\,\,2}$&$\widehat{\rho}_{2}^{\,\,2}$&  Truncation
  & $\widehat{\rho}_{1}^{\,\,2}$&$\widehat{\rho}_{2}^{\,\,2}$ &Truncation &$\widehat{\rho}_{1}^{\,\,2}$&$\widehat{\rho}_{2}^{\,\,2}$\\\hline\hline
  20& 0.5976&    0.5976  &120&    0.8883&    0.8882     &220  &0.9510&    0.8709\\
  40& 0.7446 &   0.7446&140   &0.9188 &   0.8728&240  &0.9548 &   0.8711\\
  60  &0.8107  &  0.8107&160  &0.9315  &  0.8708&260  &0.9578  &  0.8714\\
  80  &0.8481   & 0.8480&180  &0.9400   & 0.8704&280  &0.9603   & 0.8717\\
  100 &0.8719    &0.8718&200& 0.9462    &0.8706&300   &0.9625    &0.8719\\
 \hline\hline
\end{tabular}
\end{scriptsize}

 \begin{scriptsize}
\begin{tabular}{||c|c|c||||c|c|c||}
  \hline
  Truncation &$\widehat{\rho}_{1}^{\,\,2}$&$\widehat{\rho}_{2}^{\,\,2}$&Truncation &$\widehat{\rho}_{1}^{\,\,2}$&$\widehat{\rho}_{2}^{\,\,2}$
    \\\hline\hline
320 &0.9642&    0.8721&480& 0.9721   & 0.8729\\
340  &0.9658 &   0.8724&500&    0.9726&    0.8729\\
360 &0.9671  &  0.8725&520& 0.9731 &   0.8729\\
380 &0.9682   & 0.8726&540& 0.9735  &  0.8729\\
400 &0.9692    &0.8727&560& 0.9738   & 0.8729\\
420 &0.9701&    0.8727&580& 0.9741    &0.8729\\
440 &0.9709 &   0.8728&600& 0.9743&    0.8729\\
460 &0.9715  &  0.8728&&&\\
 \hline\hline
\end{tabular}
\end{scriptsize}

After a truncation of zonal polynomials of order greater than 600 the algorithm fails and no
computational chance of reaching the estimation based on the exact polynomial distribution, is
possible.

\begin{figure}[!h]
  \begin{center}
  \includegraphics[width=13cm,height=8cm]{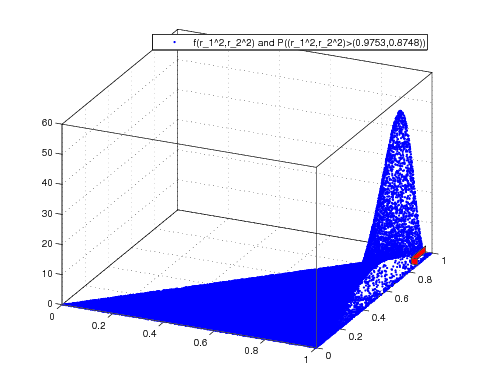}
  \caption{Eulerian polynomial joint distribution of template vs handwritten digit 3 and a
  test for a randomly selected handwritten code.}\label{fig:templatevsdigit}
  \end{center}
\end{figure}

Summarizing we have arrive at two estimates: Template vs Handwritten digit 3,
$\widehat{\rho}_{1}^{\,\,2}=0.9542$ and $\widehat{\rho}_{2}^{\,\,2}=0.6740$; meanwhile Given
Handwritten digit 3 vs Handwritten digit 3, $\widehat{\rho}_{1}^{\,\,2}=0.9753$ and
$\widehat{\rho}_{2}^{\,\,2}=0.8748$, the difference is notable, but a way for comparison can be
the addressed distribution of Theorem \ref{th:polynomial}. In this case we just compute the
probability that an extreme correlation structure as $r_{1}^{2}=0.9753$ and $r_{1}^{2}=0.8748$
due to a particular given handwritten 3, occurs under the distribution of the correlation
structure influenced by the template.  In this case, integration of distribution of Theorem
\ref{th:polynomial} provides that $P((r_{1}^{2},r_{2}^{2})>(0.9753,0.8748))<0.0015$, which
means that the Eulerian shape of the given handwritten digit 3 is so far from the usual
template of number 3. Figure \ref{fig:templatevsdigit} shows the associated distribution of the
template and an approximate region for integration.

\subsection{An application in machine vision}

\begin{figure}[!h]
  \begin{center}
  \includegraphics[width=13cm,height=8cm]{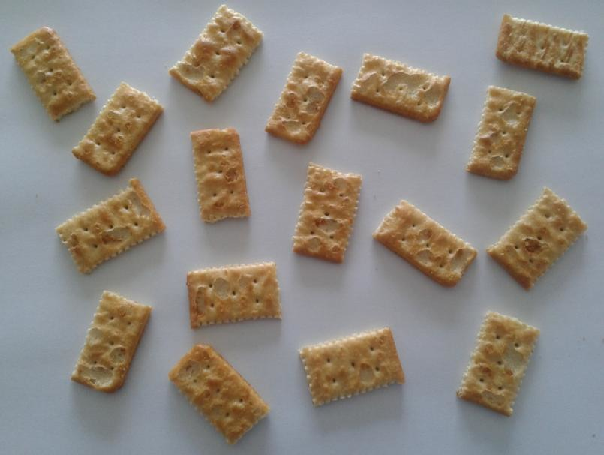}
  \caption{A 17-sample of crackers type A randomly selected, original objects.}\label{fig:sample1}
  \end{center}
\end{figure}

\begin{figure}[!h]
  \begin{center}
  \includegraphics[width=13cm,height=8cm]{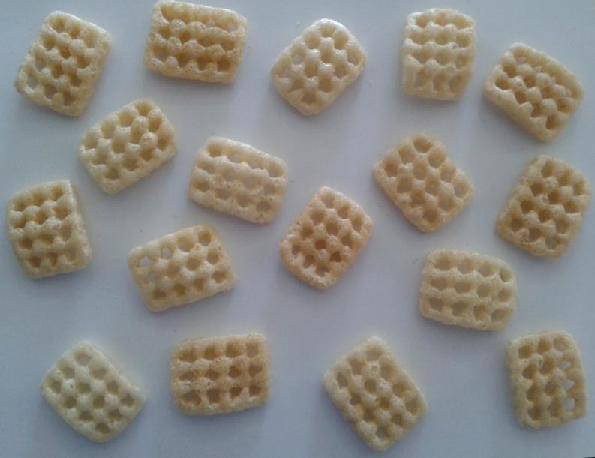}
  \caption{A 17-sample of crackers type B randomly selected, original objects.}\label{fig:sample2}
  \end{center}
\end{figure}

We end this section by studying the Eulerian shape differentiation between two populations, the
application can be useful in machinery vision. Consider two samples of 17 crackers as shown in
figures \ref{fig:sample1} (cracker type A) and \ref{fig:sample2}(cracker type B).

\begin{figure}[!h]
  \begin{center}
  \includegraphics[width=13cm,height=8cm]{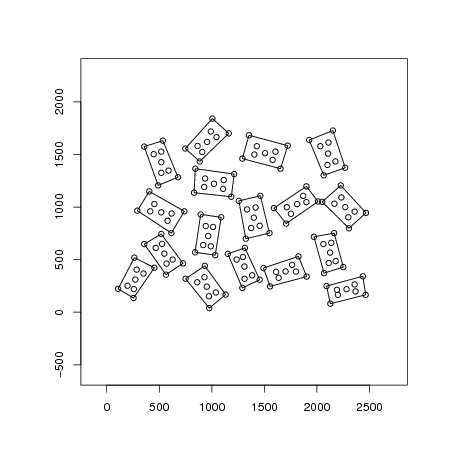}
  \caption{A 17-sample of crackers type A randomly selected with 9 landmarks.}\label{fig:landmarksA}
  \end{center}
\end{figure}

\begin{figure}[!h]
  \begin{center}
  \includegraphics[width=13cm,height=8cm]{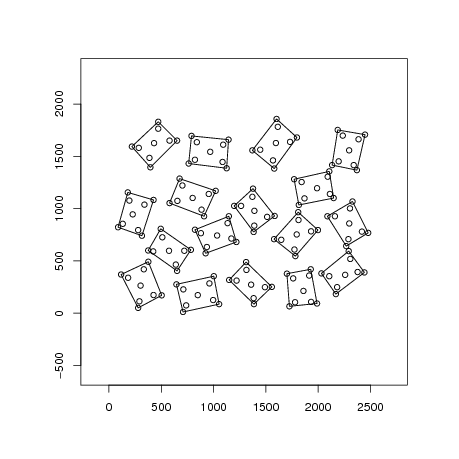}
  \caption{A 17-sample of crackers type B randomly selected with 9 landmarks.}\label{fig:landmarksB}
  \end{center}
\end{figure}

In cracker type A, we register  9 landmarks by hand, the five evident marks on the cracker and
the four corners, see Figure \ref{fig:landmarksA}. We try to place similar landmarks on cracker
type B, corresponding to the four corners, the center and the four nearest holes to the
corners, , see Figure \ref{fig:landmarksB}. It is clear the notorious variance of placing this
landmarks on cracker B versus the more easily registration on cracker A. For a complete process
in machine vision, we can also use automatic corner detection as Harris' method for the
landmark registration, but we want to increase intentional the difficulty of the pursuit
population discrimination. Note also the adequacy of objects for implementing the Eulerian
differentiation, the figures have strong symmetry and some flexibility of choosing the landmark
labels is allowed because the geometric invariants of such shape; meanwhile, in affine or
similarity studies, the labeling have more radical constraint.

We ask for a significant difference in Eulerian shape between crackers A and B. Next we list
the results of some approaches.

\textit{METHOD 1: USING A DEGENERATED TEMPLATE:}

\begin{figure}[!h]
  \begin{center}
  \includegraphics[width=13cm,height=8cm]{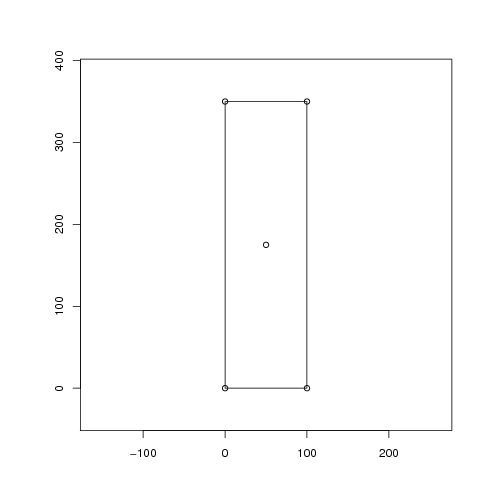}
  \caption{A degenerated template}\label{fig:template2}
  \end{center}
\end{figure}

\begin{figure}[!h]
  \begin{center}
  \includegraphics[width=13cm,height=8cm]{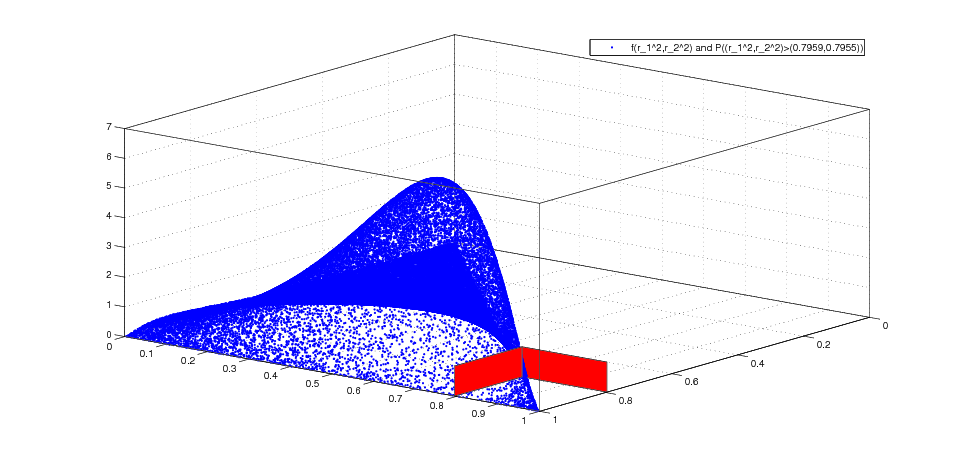}
  \caption{Eulerian polynomial joint distribution of a degenerated template vs cracker B and a test for a cracker A.}\label{fig:template2AndB}
  \end{center}
\end{figure}

We proceed with a similar methodology of the handwritten digit 3 landmark data. In this case
the selected template is just a rectangle with its center, in such way that 5 of the 9
landmarks are concentrated in the geometrical center and the remaining ones are the vertices,
see Figure \ref{fig:template2}. Corresponding estimates for population of cracker A are
$\widehat{\rho}_{1}^{\,\,2}=0.7959$ and $\widehat{\rho}_{2}^{\,\,2}=0.7955$; meanwhile
estimates for population of cracker B are, $\widehat{\rho}_{1}^{\,\,2}=0.6218$ and
$\widehat{\rho}_{2}^{\,\,2}=0.6212$, as we expect, the template is a little more correlated
with cracker A than cracker B. In this case we just compute the probability that an extreme
correlation structure as $r_{1}^{2}=0.7959$ and $r_{1}^{2}=0.7955$ due to population
correlation structure estimates of cracker B with the given template, occurs under the
distribution of the correlation structure influenced by the template and population of cracker
A. After  integration of distribution of Theorem \ref{th:polynomial} we have that
$P((r_{1}^{2},r_{2}^{2})>(0.7959,0.7955))< 0.0141$, which means that the Eulerian shape of
cracker B is different from the Eulerian shape of cracker A, even under such degenerated
template with more than 50 percent of coincident landmarks. Figure \ref{fig:template2AndB}
shows the associated distribution of the template and an approximate region for integration.

\textit{METHOD 2: USING A SIMILAR TEMPLATE TO CRACKER A:}

\begin{figure}[!h]
  \begin{center}
  \includegraphics[width=13cm,height=8cm]{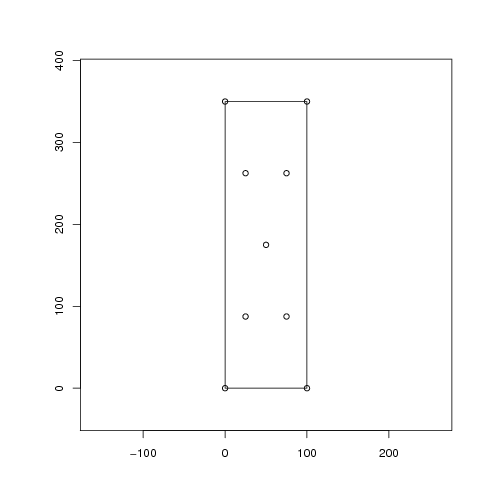}
  \caption{A template close to cracker A.}\label{fig:template1}
  \end{center}
\end{figure}

\begin{figure}[!h]
  \begin{center}
  \includegraphics[width=13cm,height=8cm]{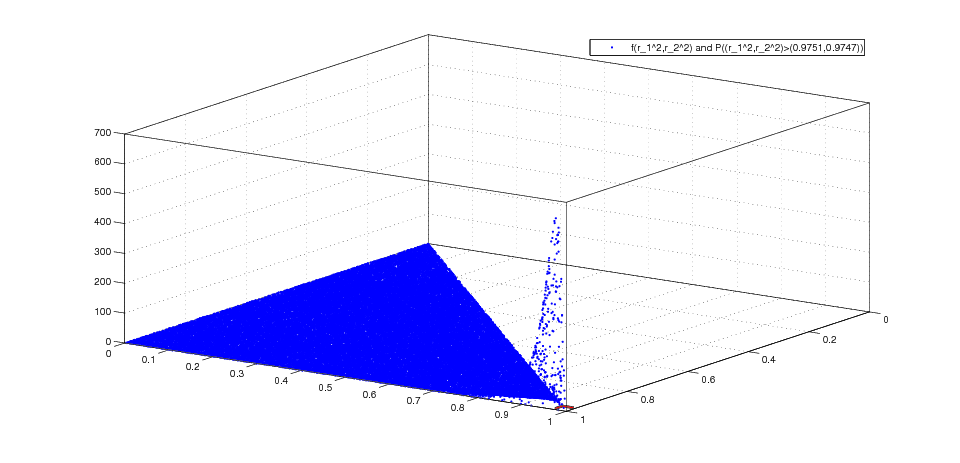}
  \caption{Eulerian polynomial joint distribution of a cracker B and template similar to cracker A,
  and a test for  cracker A.}\label{fig:template1AndB}
  \end{center}
\end{figure}

If we select a template where 4 of the above coincident landmarks are equally spaced in the
diagonals of the rectangle (very similar to cracker A), see Figure \ref{fig:template1}, the
corresponding estimates for population of cracker A are $\widehat{\rho}_{1}^{\,\,2}=0.9751$ and
$\widehat{\rho}_{2}^{\,\,2}=0.9747$; meanwhile estimates for population of cracker B are,
$\widehat{\rho}_{1}^{\,\,2}=0.9714$ and $\widehat{\rho}_{2}^{\,\,2}=0.9710$. Then we compute
the probability that an extreme correlation structure as $r_{1}^{2}=0.9751$ and
$r_{1}^{2}=0.9747$ due to population correlation structure estimates of cracker A with the
given template, occurs under the distribution of the correlation structure influenced by the
template and population of cracker B.  After  integration of distribution of Theorem
\ref{th:polynomial} we have that $P((r_{1}^{2},r_{2}^{2})>(0.9751,0.9747))<0.0017$, which means
that the Eulerian shape of cracker A is so far from the Eulerian shape of cracker B. Figure
\ref{fig:template1AndB} shows the associated distribution of the template and an approximate
region for integration.

\textit{METHOD 3: USING A TEMPLATE CLOSE TO CRACKER B:}

\begin{figure}[!h]
  \begin{center}
  \includegraphics[width=13cm,height=8cm]{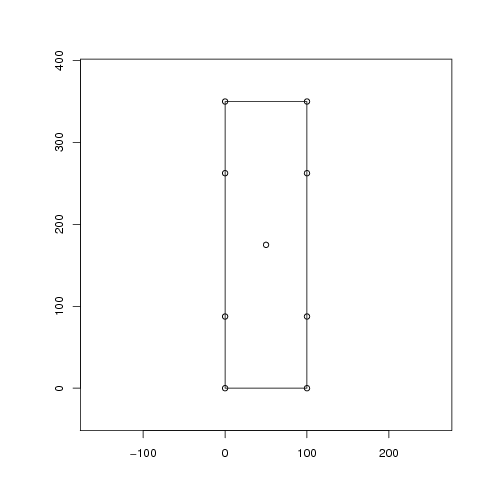}
  \caption{A template close to cracker B.}\label{fig:template3}
  \end{center}
\end{figure}

\begin{figure}[!h]
  \begin{center}
  \includegraphics[width=13cm,height=8cm]{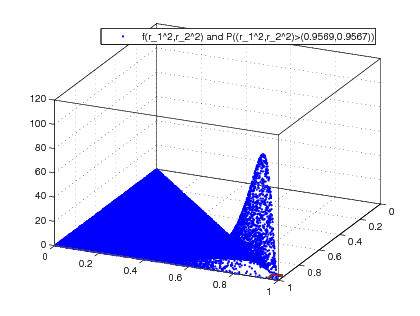}
  \caption{Eulerian polynomial joint distribution of a cracker A and a template similar to cracker B,
  and a test for  cracker B.}\label{fig:template3AndB}
  \end{center}
\end{figure}

If we select a template where 4 of the 5 coincident landmarks of method 1, are place on the
sides of the rectangle (more closer to cracker B), see Figure \ref{fig:template3} the
corresponding estimates for population of cracker A are $\widehat{\rho}_{1}^{\,\,2}=0.8928$ and
$\widehat{\rho}_{2}^{\,\,2}=0.8927$; meanwhile estimates for population of cracker A are,
$\widehat{\rho}_{1}^{\,\,2}=0.9569$ and $\widehat{\rho}_{2}^{\,\,2}=0.0.9567$. Then we compute
the probability that an extreme correlation structure as $r_{1}^{2}=0.9569$ and
$r_{1}^{2}=0.9567$ due to population correlation structure estimates of cracker b with the
given template, occurs under the distribution of the correlation structure influenced by the
template and population of cracker A.  After  integration of distribution of Theorem
\ref{th:polynomial} we have that $P((r_{1}^{2},r_{2}^{2})>(0.9569,0.9567))<0.00038$, which
means that the Eulerian shape of cracker A is so far from the Eulerian shape of cracker B.
Figure \ref{fig:template3AndB} shows the associated distribution of the template and an
approximate region for integration.

Note that the three methods ratifies the same conclusion about Eulerian shape difference, which
seems to be independent on the template (of course it must satisfies the geometrical invariants
demanding for Eulerian shape). In some experiments is difficult to suggest a template, but the
addressed Heuristic independence shows that it is not a matter of big difference among the two
pairs of estimates, but the joint distribution tends to be so narrow when the estimates are
close as in method 2, meanwhile a degenerated  template will provide apart enough the pairs of
estimates but the joint distribution for testing is so wide, as we can see in method 3.  Medium
difference as in method 1 has a proportional joint distribution size. The most important fact
is that in presence of shape differences, the selected template emerges as a theoretical or
artificial object to perform the test, but it is irrelevant (in the sense of Eulerian
restrictions). A very far artificial template will separate differences but the associated
distribution has big variance, and it affects the resolution of the test. Some mixing ideas for
an artificial template close to intuition could be a good choice. Another usual conjecture for
the template is a real observation of the sample, however, we have to be sure that the assigned
template really belongs to the addressed population.

We try to apply this methodology with a very sensitive landmark data, the schizophrenia
landmark date due to \citet{B:96}. In this case 13 landmarks in 2 dimensions, were register on
a near midsaggital brain scan of 14 controls and  14 patients with schizophrenia. The landmarks
choosen in the  MRI  were: (1) splenium, posteriormost point on corpus callosum; (2) genu,
anteriormost point on corpus callosum; (3) top of corpus callosum, uppermost point on arch of
callosum (all three to an approximate registration on the diameter of the callosum); (4) top of
head, a point relaxed from a standard landmark along the apparent margin of the dura; (5)
tentorium of cerebellum at dura; (6) top of cerebellum; (7) tip of fourth ventricle; (8) bottom
of cerebellum; (9) top of pons, anterior margin; (10) bottom of pons, anterior margin; (11)
optic chiasm; (12) frontal pole, extension of a line from landmark 1 through landmark 2 until
it intersects the dura; (13) superior colliculus; see \citet{B:96} and \citet{DM98}.

\begin{figure}[!h]
  \begin{center}
  \includegraphics[width=13cm,height=8cm]{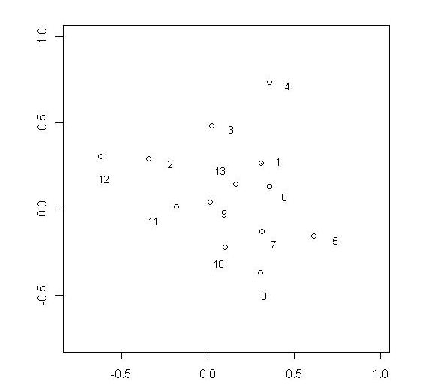}
  \caption{A template corresponding to a randomly sampled schizophrenic patient.}\label{fig:templatesquizo}
  \end{center}
\end{figure}

First we use as a template one of the MRI sampled from schizophrenic patients, see Figure
\ref{fig:templatesquizo}. The computations are simpler and no new integrals must be computed in
Theorem \ref{th:polynomial}, because this data base has the same number of landmarks than the
previously studied problem of handwritten digit 3. Again the distribution of the squared
canonical correlation population is just a polynomial of degree 12 in the latent roots of the
matrix argument.

\begin{figure}[!h]
  \begin{center}
  \includegraphics[width=13cm,height=8cm]{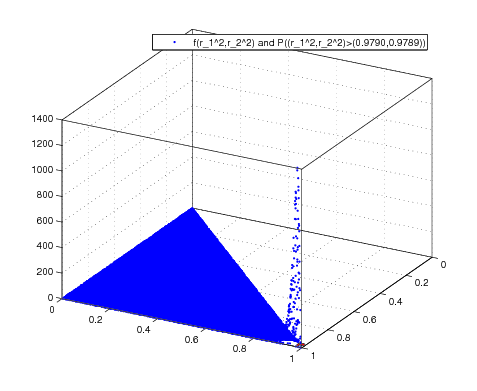}
  \caption{Eulerian polynomial joint distribution of controls and a schizophrenic patient-template,
  and a test for schizophrenic patients.}\label{fig:templatesquizovscontrols}
  \end{center}
\end{figure}

The corresponding estimates for population of schizophrenic patients are
$\widehat{\rho}_{1}^{\,\,2}=0.9790$ and $\widehat{\rho}_{2}^{\,\,2}=0.9789$; meanwhile
estimates for population of controls are a little less as we expect,
$\widehat{\rho}_{1}^{\,\,2}=0.9783$ and $\widehat{\rho}_{2}^{\,\,2}=0.9781$. Then we compute
the probability that an extreme correlation structure as $r_{1}^{2}=0.9790$ and
$r_{1}^{2}=0.9789$ due to population correlation structure estimates of schizophrenic patients
with the given template, occurs under the distribution of the correlation structure influenced
by the template and population of controls.  After  integration of distribution of Theorem
\ref{th:polynomial} we have that $P((r_{1}^{2},r_{2}^{2})>(0.9790,0.9789))<0.0024$, which means
that the Eulerian shape of schizophrenic patients is different from the Eulerian shape of
controls. Figure \ref{fig:templatesquizovscontrols} shows the associated distribution of the
template and an approximate region for integration.

A similar result is obtained by taking an arbitrary template randomly selected from the sample of controls.

\textit{Concluding remark:} Finally, observe that the exact Eulerian polynomial joint
distribution can be used for testing the structure of correlation, as in the same way that
equality mean shape tests are usually done in shape analysis, with the difference that we are
using the exact distribution instead of the asymptotic distributions of the statistics
traditional implemented. Of course, the inference on correlation structure differs from the
inference about mean shape, and both analysis complement each other, but it seems that the
Eulerian shape analysis is sufficient for detecting differences between populations as the
similar mean shape analysis. In fact, we have applied the methods presented here to some of the
classical landmark data studied in literature by asymptotic distributions or inference in the
tangent plane as those presented in \citet{DM98}, and we have obtained the same conclusions
about mean shape differences, among the applications we do not include in this manuscript we
count the following landmark data: mouse vertebra, gorilla skulls and macaque skulls. However,
for completeness, a subsequent work will include the addressed associated Eulerian mean shape
theory and we will back to this applications.

Observe, that the exact Eulerian polynomial joint distribution can be used for testing the
structure of correlation, as in the same way that equality mean shape tests are usually done in
shape analysis, with the difference that we are using the exact distribution instead of the
asymptotic distributions of the statistics traditional implemented. Of course, the inference on
correlation structure differs from the inference about mean shape, and both analysis complement
each other.  In fact, a subsequent work will include the addressed associated Eulerian mean
shape theory and we will back to this applications.

\section*{Acknowledgments}
F. Caro was supported by University of Medellín in the context of a join research project with
University of Toulouse and University of Bordeaux.

\end{document}